\documentclass{amsart}
\usepackage{amsmath,latexsym,amssymb,amsfonts}
\usepackage{amscd,amssymb,epsfig, epsf}

\newtheorem{theorem}{Theorem}
\newtheorem{lemma}[theorem]{Lemma}

\newtheorem{MainThm}{Theorem}

\title{Free Infinite Divisibility for
Q-Gaussians}
\author{M. Anshelevich}
\address{Department of Mathematics,
Texas A\&M University, 
College Station TX 77843-3368, U.S.A.}
\email{manshel@math.tamu.edu}
\author{S. T. Belinschi}
\address{Department of Mathematics \& Statistics,
University of Saskatchewan, and Institute of Mathematics ``Simion
Stoilow'' of the Romanian Academy;
106 Wiggins Road
Saskatoon, SK S7N 5E6, CANADA} 
\email{belinschi@math.usask.ca}
\author{M. Bo\.zejko}
\address{Instytut Matematyczny, Uniwersytet Wroc{\l}awski, Pl. Grunwaldzki 2/4,50-384, Wroc{\l}aw, POLAND} 
\email{bozejko@math.uni.wroc.pl}
\author{F. Lehner}
\address{Institute for Mathematical Structure Theory, Graz Technical 
University, Steyrergasse 30, 8010 Graz, AUSTRIA} 
\email{lehner@finanz.math.tu-graz.ac.at}

\date{\today}
\begin{document}

\thanks{Work of M. Anshelevich was supported in part by NSF grant DMS-0900935.}
\thanks{Work of S. Belinschi was supported by a Discovery
Grant from NSERC and a University of Saskatchewan start-up grant.}
\thanks{Work of M. Bo\.zejko was partially supported by Polish Ministry of
Science and Higher Education by grant no. N N201 364436.}

\begin{abstract}We prove that the $q$-Gaussian distribution is freely
infinitely divisible for all $q\in[0,1]$.
\end{abstract}

\maketitle

\section{Introduction}

In this note we prove that the $q$-Gaussian distribution
introduced by Bo\.zejko and Speicher in \cite{BS} (see also the paper \cite{BKS} of
Bo\.zejko, K\"{u}mmerer and Speicher) is freely infinitely
divisible when $q\in[0,1]$.

We shall give a short outline for the context of this problem.
In probability theory, the class of infinitely divisible distributions plays a crucial role, in the study of limit theorems, L{\'e}vy processes etc. So it was a remarkable discovery of Bercovici and Voiculescu
\cite{BVIUMJ} that there exists a corresponding class of \emph{freely} infinitely divisible distributions in 
free probability. These distributions are typically quite different from the classically infinitely divisible ones; for example, many of them are compactly supported. Nevertheless, work of numerous authors culminating in the paper by \cite{BP} 
Bercovici and Pata showed that free ID distributions are in a precise bijection with
the classical ones, this bijection moreover having numerous strong properties. For example, the semicircular law is the free analog of the normal distribution. From this bijection, one might get the intuition that, perhaps with very rare exceptions such as the Cauchy and Dirac distributions, 
some measures belong to the
``classical'' world and some to the ``free'' world. However, \cite[Corollary 3.9]{BBLS} indicates that this
intuition may be misleading: the normal distribution, perhaps the most important among the classical 
ones, is also freely infinitely divisible.

One approach towards understanding the relationship between the classical and free probability theories, and in particular the Bercovici-Pata bijection, have been attempts to construct an interpolation between these two theories. The oldest such construction, due to Bo{\.z}ejko and Speicher, is the 
construction of
the $q$-Brownian motion. In particular, it provides a probabilistic interpretation for a (very classical)
family of $q$-Gaussian distributions, which interpolate between the normal ($q=1$) and the semicircle ($q=0$) laws. Probably the best known description of the $q$-Gaussian distributions is in terms of 
their orthogonal polynomials $H_n(x|q)$ - called the $q$-Hermite polynomials - defined by the 3-term
recurrence relation $xH_n(x|q)=H_{n+1}(x|q)+\frac{1-q^n}{1-q}H_{n-1}(x|q)$, with initial conditions
$H_0(x|q)=1,$ $H_1(x|q)=x$. These polynomials have been studied for a long time: probably their first 
appearance under this guise occurred in L. J. Rogers' 1893 paper \cite{LGR}. 
Bo\.zejko and Speicher  construct their example of $q$-Brownian motion using 
creation-annihilation operators on a ``twisted'' Fock space: given a separable Hilbert space $H$ and 
$f\in H$, we let $c^*(f)$ be the
left creation operator on the Fock space $\mathcal F(H)$ and $c(f)$ be its adjoint. 
The authors provide in \cite{BS} a scalar product $(\ \cdot\ |\ \cdot\ )_q$ 
on (a quotient space of) $\mathcal F(H)$ for which the so-called $q$-canonical 
commutation relation holds:
$$
c(f)c^*(g)-qc^*(g)c(f)=(f|g)_q{\bf 1},\quad f,g\in H,
$$
where {\bf 1} is the identity operator on $\mathcal F(H)$. 
It is a fundamental result  that, when $(f|f)_q=1$, 
the distribution of the self-adjoint operator $c(f)+c^*(f)$ with 
respect to the vacuum state on $\mathcal F(H)$ is the (centered) $q$-Gaussian distribution of variance
one, having the $q$-Hermite polynomials as orthogonal polynomials (see \cite[Theorem 1.10]{BKS}). 
In this case, as mentioned before, when $q=1$, $c(f)+c^*(f)$ is distributed according to the classical 
normal distribution $(2\pi)^{-1/2}\exp(t^2/2)dt$, and when $q=0$ according to the free central 
limit - the Wigner distribution with density $(2\pi)^{-1}\sqrt{4-t^2}\chi_{[-2,2]}(t)$. A formula for the 
density $f_q$ of the $q$-Gaussian distribution is provided in \cite{LM}:
$f_q(x)=\pi^{-1}\sqrt{1-q}\sin\theta\prod_{n=1}^\infty(1-q)^n|1-q^ne^{2i\theta}|^2,$ where
$\cos\theta=\frac{x}{2}\sqrt{1-q}.$ The interested reader might want to note \cite[Section III]{LM} that,
up to a multiplicative constant depending only on $q$, $f_q$ is a theta function.
For numerous details on properties and applications of $q$-Gaussian processes we refer to \cite{BKS}
and references therein.

As the construction described above suggests, $q$-Gaussians provide useful examples
in operator algebras. The  von Neumann algebras generated by families of
$q$-Gaussians are shown  to exhibit several  interesting properties, and we shall
list a few below; however, the structure of these algebras still remains largely mysterious. 
It is shown in \cite{Nou} that algebras generated by such families when $f$ runs through $H$ and 
dim$(H)\ge2$ are non-injective; the paper \cite{ER}
proves that algebras generated by at least two $q$-Gaussians corresponding to orthogonal $f$s 
are always factors when $|q|<1$ (see also \cite{H,BKS,Sniany}), and in \cite{Dima} Shlyakhtenko
provides estimates for the non-microstates free entropy of $n$-tuples of such $q$-Gaussians, 
estimates which guarantee that the algebra they generate is solid in the sense of Ozawa whenever
$q<\sqrt2-1$. Moreover, recently Bo\.zejko \cite{B} proved that in von Neumann algebras generated by two
$q$-Gaussians, the Bessis-Moussa-Villani conjecture holds.

There are several strictly probabilistic approaches to $q$-Gaussians: we would like to mention
the stochastic integration with respect to $q$-Brownian motion of Donati-Martin
\cite{CDM}, the $q$-deformed cumulants, a $q$-convolution defined on a restricted class of
probability measures and
$q$-Poisson processes studied in \cite{A}
and a random matrix model provided in \cite{Kemp}. 

However, classical or free probability aspects of $q$-Gaussians have been less studied.
In this paper, we show that all of these distributions, for $0 \leq
q \leq 1$, are freely infinitely divisible. This may be an indication that the class of freely infinitely divisible distributions, despite the Bercovici-Pata bijection, is quite different from the classical one, and is yet to be understood completely.
The conjecture that $q$-Gaussian distributions are freely infinitely divisible when $q\in[0,1]$, formulated 
initially by the two last named authors and R. Speicher, was motivated among others by the recently proved free infinite divisibility
of the classical Gaussian (corresponding to $q=1$).
It has been shown in \cite{BBLS} that the Askey-Wimp-Kerov
distributions with parameter $c\in[-1,0]$, and in particular
the classical normal distribution (corresponding to $c=0)$
are $\boxplus$-infinitely divisible. This provided free infinite
divisibility for distributions of several noncommutative Brownian
motions (see for example \cite{BS2,BDJ,Buchholz}), interpolating between the 
classical central limit (the
normal distribution) and the free central limit (the Wigner semicircle
law). However, until now it remained a mystery whether this first
(and most famous) such example of interpolation consists also of
$\boxplus$-infinitely divisible distributions. Several numerical
verifications performed by one of us seemed to indicate
this is indeed the case. 
Here we shall give an answer to this question:
\begin{theorem}\label{main}
The $q$-Gaussian distribution $f_q(x)dx$
is freely infinitely divisible for all $q\in[0,1]$.
\end{theorem}
Our method to prove this result will be the same as in \cite{BBLS}, 
namely we will construct an inverse to the Cauchy transform $G_{f_q}$
of the $q$-Gaussian defined on the whole lower half-plane.
Then the Voiculescu transform $\phi_{f_q}(1/z)=G_{f_q}^{-1}(z)-(1/z)$
clearly has an extension to the whole complex upper half-plane 
$\mathbb C^+$.
An application of the following theorem of Bercovici and Voiculescu from \cite{BVIUMJ}
yields the desired result:

\begin{theorem}\label{bvid}
A Borel probability measure $\mu$ on the real line is $\boxplus$-infinitely divisible if and only if
its Voiculescu transform $\phi_\mu(z)$ extends to an analytic function $\phi_\mu\colon
\mathbb C^+\to\mathbb C^-$.
\end{theorem}

The paper is organized as follows: in the second section, we introduce
several notions and preliminary results used in our proof,
and in the third section we give the proof of Theorem \ref{main}.

\noindent{\bf Acknowledgements.} We would like to thank Roland Speicher for many useful discussions
regarding this paper.

\section{Preliminary results: the importance of an entire function}

It is shown in the paper \cite{Pawel} of Pawel Szab{\l}owski (see also 
\cite{MSV}) that the density of the $q$-Gaussian distribution
of mean zero and variance one,
with respect to the Lebesgue measure is given by the formula
$$
f_q(x)=\frac{\sqrt{1-q}}{2\pi}\sqrt{4-(1-q)x^2}\sum_{k=1}^\infty
(-1)^{k-1}q^{\frac{k(k-1)}{2}}\mathcal U_{2k-2}\left(\frac{x\sqrt{
1-q}}{2}\right),
$$
for $|x|\le2/\sqrt{1-q}$, where $\mathcal U_k$ is the $k^{\rm th}$ Chebyshev polynomial of the second kind - defined by the relation $\mathcal U_k(\cos\theta)=\frac{\sin((k+1)\theta)}{\sin\theta}$. 
In our present paper we will consider 
this as being the definition of the $q$-Gaussian distribution.
For simplicity of notation, we will re-normalize this density
to being supported on $[-2,2]$, by replacing $x$ with $x/\sqrt{1-q}:$
$$
f_q(x)=\frac{1}{2\pi}\sqrt{4-x^2}
\sum_{k=1}^\infty
(-1)^{k-1}q^{\frac{k(k-1)}{2}}\mathcal U_{2k-2}\left(\frac{x}{2}
\right),\quad -2\leq x\leq 2.
$$

Recall that the Cauchy (or Cauchy-Stieltjes) transform of a Borel
probability measure $\mu$ on $\mathbb R$ is by definition
$$
G_\mu(z)=\int_\mathbb R\frac{1}{z-x}\,d\mu(x),\quad z\in
\mathbb C\setminus\mathbb R.
$$
This function maps the upper half-plane $\mathbb C^+$ into the 
lower half-plane $\mathbb C^-$, satisfies $G_\mu(\overline{z})
=\overline{G_\mu(z)}$, and extends analytically through
the complement of the support of $\mu$. Also, of some importance for us will be 
the map $F_\mu(z)=\frac{1}{
G_\mu(z)}$. This function satisfies the inequality $\Im F_\mu(z)>\Im z$, $z\in\mathbb C^+$,
whenever $\mu$ is not a point mass. For more details
we refer the reader to the third chapter of \cite[Chapter III]{akhieser}.

Integrating $(z-x)^{-1}f_q(x)$ with respect to the Lebesgue 
measure on $[-2,2]$, we obtain for any $z\in\mathbb C^+$ that 
\begin{eqnarray*}
G_{f_q}(z) & = & \int_{-2}^2\frac{1}{z-x}f_q(x)dx\\
& = & \sum_{k=1}^\infty
(-1)^{k-1}q^{\frac{k(k-1)}{2}}\left(\frac{1}{2\pi}
\int_{-2}^2\frac{1}{z-x}\mathcal U_{2k-2}\left(\frac{x}{2}\right)
\sqrt{4-x^2}\,dx\right)\\
& = & 
\sum_{k=1}^\infty
(-1)^{k-1}q^{\frac{k(k-1)}{2}}G_s(z)^{2k-1}\\
& =  & \sum_{k=0}^\infty
(-1)^{k}q^{\frac{k(k+1)}{2}}G_s(z)^{2k+1}
,
\end{eqnarray*}
where $G_s$ is the Cauchy transform of the semicircular law:
$$
G_s(z)=\frac{1}{2\pi}\int_{-2}^2\frac{1}{z-x}\sqrt{4-x^2}\,dx
=\frac{z-\sqrt{z^2-4}}{2},
\quad z\in\mathbb C^+.
$$
This function has an analytic extension to the lower half-plane
$\mathbb C^-$ through the interval $(-2,2)$ that does {\em not} coincide with
$(2\pi)^{-1}\int_{-2}^2(z-x)^{-1}\sqrt{4-x^2}dx$ when $z\in\mathbb C^-$:
when we consider
the same branch of the square root as above, it is of the form
$G_s(z)=\frac{z+\sqrt{z^2-4}}{2}$ (meaning, the asymptotics of this
extension at $-i\infty$ is $O(z)$). An analysis of these
two formulas guarantee us that $G_s$ maps $\mathbb C^+\cup(-2,2)\cup
\mathbb C^-$ bijectively onto $\mathbb C^-$ by mapping $\mathbb C^+$
into the lower half of the unit disc $\mathbb D$ (the piece
of $\partial G_s(\mathbb C^+)$ that forms the interval $[-1,1]$
is $G_s([-\infty,-2]\cup[2,+\infty])$, with the two infinities
identified and $G_s(\infty)=0$), while $G_s(\mathbb C^-)$ is the 
complement of $\overline{G_s(\mathbb C^+)}$ in $\mathbb C^-$.
In addition, $G_s(i\mathbb R_+)=i[-1,0)$, $G_s(i\mathbb R_-)=i
(-\infty,-1]$ and, when discarding the $i$, $G_s$ is monotonic
increasing on the imaginary axis.
These simple observations will be essential in our proof.

Now observe that the above remarks translate into 
$G_{f_q}=g_q\circ G_s$, where 
\begin{equation}\label{1}
g_q(w)=\sum_{k=0}^\infty
(-1)^{k}q^{\frac{k(k+1)}{2}}w^{2k+1},\quad w\in\mathbb C,
\end{equation}
is an entire function for any $q\in[0,1)$ (in fact for any $|q|<1$). \
(The reader will note that in terms of basic hypergeometric
functions $g_q$ can be written as 
$g_q(w)=w\cdot{}_1\phi_0(q^\frac12|q^\frac12,(iq^\frac12w)^2)$; however, we shall not use
this fact directly in our proof.)
We list below a few properties  of $g_q$ which will be used in our proof:
\begin{enumerate}
\item $\lim_{q\to0} g_q(w)=w,w\in\mathbb C$, and the limit is uniform 
on compacts in $\mathbb C$;
\item $\lim_{q\to1} g_q(w)=\frac{w}{w^2+1},w\neq\pm i$;
\item $g_q(-\overline{w})=-\overline{g_q(w)};$
\item $g_q(\overline{w})=\overline{g_q(w)}$ - in particular
$g_q(\mathbb R)\subseteq\mathbb R$;
\item $g_q(ic)=i\sum_{k=0}^\infty q^\frac{k(k+1)}{2}c^{2k+1}$
for any real $c$, so that $g_q(i(-\infty,0])=i(-\infty,0]$,
$g_q(i[0,+\infty))=i[0,+\infty)$, with $g_q(i(\pm\infty))=
i(\pm\infty)$, $g_q(0)=0$, and $g_q$ is monotonic on the imaginary 
axis (in the same sense as $G_s$ is);
\item $g_q'(0)=1$.
\end{enumerate}
The limits in items (1), (2) above being uniform on 
the corresponding compact subsets, we can conclude that, given $M>0$, for $q>0$ sufficiently small, 
depending on $M$, we have that $g_q$ is injective on the ball of radius $M$ centered
at the origin, In addition, Equation \eqref{real} below guarantees that $g_q$ has no limit at infinity
along the real axis; items (3) and (4) above together with this remark guarantee that 
$g_q'$ will have at least two zeros on the real
line, symmetric with respect to zero, if $q>0$. 

It is time for stating a few results to be used in our proof:

\begin{theorem}\label{order}
A necessary and sufficient condition that $$f(z)=\sum_{n=0}^\infty
a_nz^n$$ should be an entire function of order $\varrho$ is that
$$
\liminf_{n\to\infty}\left(\frac{-\log|a_n|}{n\log n}\right)=
\frac1\varrho.
$$
\end{theorem}
(This is \cite[Theorem 4.12.1]{Holland}.)
Thus, for our $g_q$, we can take $\infty=\lim_{n\to\infty}
\left(\frac{-n(n+1)\log q}{2(2n+1)\log(2n+1)}\right)=
\frac1\varrho$ to conclude that the order $\varrho$ of $g_q$ is 
zero.

Let us also give the Denjoy-Carleman-Ahlfors Theorem:
\begin{theorem}\label{DCA}
\begin{enumerate}
\item[{\rm(i)}] An entire function of order $\varrho$ has at most $2
\varrho$ finite asymptotic values.
\item[{\rm(ii)}] For an entire function of order $\varrho$ the sets
$\{z\in\mathbb C\colon|f(z)|>c\}$, $c\ge0$, have at most $\max\{2
\varrho,1\}$ components.
\end{enumerate}
\end{theorem}
In particular, the function $g_q$ has no finite asymptotic value and the sets
$\{z:|g_q(z)|>c\}$ have exactly one connected component.

Two very famous theorems, available in any complex analysis text:
\begin{theorem}\label{open}
Any nonconstant analytic function is open, i.e. maps open sets into
open sets. In particular, if $f$ is analytic and nonconstant, 
$\partial f(D)\subseteq f(\partial D)$ for a given open $D$.
\end{theorem}

\begin{theorem}\label{var}
Assume that $f,g$ are analytic on the simply connected domain
$D$, and $\gamma$ is a rectifiable simple closed curve in $D$. If
$|f(z)+g(z)|<|f(z)|+|g(z)|$ for all $z$ in the range of $\gamma$,
then $f$ and $g$ have the same number of zeros 
in the subdomain delimited by $\gamma$ inside $D$.
\end{theorem}
The first is the open mapping theorem, the second a weakened
version of Rouch\'e's theorem - the variation of argument.

Finally, let us denote following \cite{AE}, page 81,
\begin{equation}\label{theta}
\theta_1(z)=-iq^{1/4}\sum_{n=-\infty}^\infty(-1)^n q^{n^2+n}
e^{i(2n+1)z}=-iqG\prod_{n=1}^\infty(1-q^{2n-2}e^{-2iz})(1-q^{2n}
e^{2iz})e^{iz}.
\end{equation}
We do the obvious substitution $z=-i\log w$ to get the new function
\begin{eqnarray}
\Theta(w) & = & -iq^{1/4}\sum_{n=-\infty}^\infty (-1)^n
q^{n^2+n}w^{2n+1}\nonumber\\
& = & -iq^{1/4}\left(g_{q^{2}}(w)-g_{q^{2}}\left(\frac1w\right)
\right).
\end{eqnarray}
The letter
$G$ in formula \eqref{theta} denotes a complex constant, not a Cauchy transform; 
we preserve this notation just in order to follow \cite{AE}. Second, while $\theta_1$ is
entire (and one can see that from either of the two formulas - 
the infinite product or the bi-infinite sum - since, for example,
$|(e^{iz})^{(2\pm\frac1n)}|<|q^{-|n|+1}|$ for $|n|$ large enough),
this is not the case, of course, with $\Theta$.
(We will generally ignore reparametrizations of $q$, since they 
won't be important.) This new function is not analytic anymore at
zero. Here it is important however to observe the product formula
for our function: the above provides us with
\begin{equation}\label{real}
\Theta(w)=-iq^{1/4}\left(g_{q^{2}}(w)-g_{q^{2}}\left(\frac1w\right)
\right)
=-iqG\prod_{n=1}^\infty(1-q^{2n-2}w^{-2})(1-q^{2n}
w^{2})w.
\end{equation}
Thus, all zeros of $\Theta:w\mapsto -iq^{1/4}\left(g_{q^{2}}(w)-g_{q^{2}}
\left(\frac1w\right)\right)$ are real (namely $w=\pm q^{n-1}$
and $w=\pm q^{-n}$, $n\in\mathbb N$; $w=0$ is not a zero of this
function, since
zero is an accumulation point of other zeros).

Let us give a few trivial lemmas:
\begin{lemma}\label{curve}
Assume that an entire analytic function $h$ has order zero. Then the 
preimage $h^{-1}(l)$ of any piecewise smooth curve $l$ with both ends 
at infinity has all its ends tend to infinity and for any component
$\gamma$ of $h^{-1}(l)$, we have $h(\gamma)=l$.
\end{lemma}
\begin{proof}
It is clear from Theorem \ref{open} and the definition of an entire
function that the preimage of any curve with both ends at infinity
cannot have an end in the complex plane.

On the other hand,
by the Denjoy-Carleman-Ahlfors Theorem (Theorem \ref{DCA}), we know
that there is no path $\gamma$ going to infinity so that
the limit of $h$ along $\gamma$ is finite. If we assume that
there is a branch $\gamma$ of $h^{-1}(l)$ so that $h(\gamma)\neq l$, 
then there exists a complex number $c\in l$ which is an asymptotic
value for $h$ at infinity, contradicting Theorem \ref{DCA}. 
\end{proof}

\begin{lemma}\label{picard}
\begin{enumerate} \item[(a)]
If $\gamma$ is the boundary of a simply connected domain $D$ in 
$\mathbb C$, $\gamma$ has both ends at infinity, and $h(\gamma)=\mathbb 
R$, then either $h(D)$ is one of the domains $\mathbb C^+,\mathbb C^-,
$ or $\overline{h(D)}=\mathbb C$.
\item[(b)] If $h(\partial D)$ is a half-line $s$ in $\mathbb C$, then
${h(D)}\supseteq\mathbb C\setminus s$.
\end{enumerate}
\end{lemma}
\begin{proof}
As seen in Theorem \ref{open}, $\partial h(D)\subseteq h(\gamma)=
\mathbb R$. It is clear that if $h(D)$ is not a half-plane, then
its closure must be all of $\mathbb C$. This proves (a). The proof
of (b) is identical: $\partial h(D)\subseteq h(\partial D)=s$,
so ${h(D)}\supseteq\mathbb C\setminus s$.
\end{proof}

\begin{lemma}\label{inj}
With the notations from the previous lemma, if $\gamma$ is a 
rectifiable path, $h(D)=\mathbb C^+$, and $h$ is injective on $\gamma$,
then $h$ maps $D$ conformally onto $\mathbb C^+$.
\end{lemma}
\begin{proof}
This lemma is a direct consequence of \cite{CP}, Chapters 1 and 2.
\end{proof}

\section{Proof of free infinite divisibility for $f_q(x)dx$}

We are now ready to prove the main result. For the comfort of the 
reader, we restate our main Theorem \ref{main} below.

\begin{MainThm}

The $q$-Gaussian distribution $f_q(x)dx$
is freely infinitely divisible for all $q\in[0,1]$.
\end{MainThm}
\begin{proof}
In our proof, we will follow the outline described in the introduction.
Namely, we will find a domain $X_q$ in the lower half-plane containing
the lower half of the unit disc $\mathbb D$ with the property that
$g_q(X_q)=\mathbb C^-$ and $g_q$ is injective on $X_q$. Since
we have shown that $G_s$ maps $\mathbb C^+\cup\mathbb R$ injectively 
into the closure of the lower half of the unit disc and (when 
considering the correct extension through $(-2,2)$) $\mathbb C^-$ 
injectively into $\mathbb C^-\setminus\overline{\mathbb D}$, it will 
follow that $G_{f_q}^{-1}=G_s^{-1}\circ g_q^{-1}$ extends
to $\mathbb C^-$. Thus,
$$
F_{f_q}(\cdot)=\frac{1}{g_q(G_s(\cdot))}\colon G^{-1}_s(X_q)\mapsto \mathbb C^+
$$ 
is a bijective correspondence. (The reader should keep in mind
that $G_s^{-1}$ is the inverse of the extension of $G_s|_{\mathbb C^+}$
through $(-2,2)$, so $G^{-1}_s(X_q)\supset\mathbb C^+$.) 
Since $X_q$ is included in the lower half-plane, the choice of the extension of $G_s$
guarantee that $G^{-1}_s(X_q)\cap\mathbb R=(-2,2)$. This implies the existence of $\phi_{f_q}(z)=
F_{f_q}^{-1}(z)-z$ for all $z\in \mathbb C^+$. Recalling now that $\Im F_{f_q}(w)>\Im w$
for all $w\in\mathbb C^+$, we obtain that $\Im F_{f_q}(w)>\Im w$ for all $w\in G^{-1}_s(X_q)$,
and hence $\Im \phi_{f_q}(z)=
\Im F_{f_q}^{-1}(z)-\Im z<0$ for all $z\in\mathbb C^+$.
An application of Theorem
\ref{bvid} yields the desired conclusion.

First, remark that, 
since $g_q'(0)=1$, for any fixed $q\in(0,1)$ there exist
two constants $K_q,M_q>0$ so that $g_q$ is injective on
$K_q\mathbb D$ and $g_q(K_q\mathbb D)\supset M_q\mathbb D$.

Next, we shall construct the domain $X_q$ described in the first
paragraph of the proof. Since $g_q(z)=-g_q(-z)$ and $g_q(\overline{z})
=\overline{g_q(z)}$, $z\in\mathbb C$, it will clearly be enough
to find the right side of $X_q$, as $X_q$ must be symmetric
with respect to the imaginary axis.
As suggested by the last lemma, we shall first find 
a continuous simple path inside $\overline{\mathbb C^-}$ whose image
via $g_q$ is the real line.
We observe in addition that, due to the symmetry noted in items
(3), (4) and (5) in the above list of properties of $g_q$, it will be enough to determine
the right half of such a path (see Figure \ref{fig1}.)

We start our path from zero.
Let us now ``walk'' along the real axis, in the positive direction
(by item (3) above, it is enough if we cover the right half)
until we encounter a zero of $g_q'$, call it $d_q$. Clearly,
$g_q([0,d_q])=[0,g_q(d_q)]$ is a bijective identification. Around
$d_q$, $g_q$ will be an $n$-to-one cover (where $n$ is the order of the
zero of $g_q(\cdot)-g_q(d_q)$), and we will choose the path ``first to 
the right'' which continues $g_q^{-1}([0,g_q(d_q)])$ beyond $d_q$, to 
$g_q^{-1}([0,+\infty))$.
Observe that this path escapes now in the lower half-plane and,
since $g_q(\mathbb R)\subseteq\mathbb R$, it remains in the lower
half-plane.
Now whenever we meet another such critical point for 
$g_q$ on this path $g_q^{-1}([0,+\infty))$, we turn again ``first to 
the right'', so that immediately to the right of this branch of 
$g_q^{-1}([0,+\infty))$ of ours, $g_q$ is injective. Call this 
branch - which is a rectifiable, piecewise analytic path - 
$\gamma_q$. Observe moreover that, by item (5), $\gamma_q$ is 
confined to the lower right quadrant of the complex plane (it cannot 
cut the imaginary axis, as the imaginary axis is mapped into itself).

Apriori it 
is not clear whether this branch of $g_q^{-1}([0,+\infty))$ is in 
fact existing, i.e. $g_q(\gamma_q)=[0,+\infty)$. 
However, as we have seen above, $g_q$ has order zero, so by 
Lemma \ref{curve} (the application of the 
Denjoy-Carleman-Ahlfors Theorem), $g_q$ can have no finite asymptotic 
values at infinity, so indeed $\gamma_q$ as above (i.e. $g_q(\gamma_q)=
[0,+\infty)$) exists. Thus $\gamma_q\cup(-\overline{\gamma_q})$ (here
$\overline{\gamma}$ means complex conjugate, not closure) forms the
boundary of a unique simply connected domain $X_q\subseteq\mathbb C^-$
which contains the lower half of the imaginary axis and has a
piecewise analytic, everywhere continuous, boundary. Clearly
the lower half of the unit disc is included in $X_q$. Indeed, if
$\gamma_q\cap\mathbb D\neq\varnothing$, then
there exists $z_0\in\mathbb C^+$ so that $G_s(z_0)\in
\gamma_q\cap\mathbb D$, and thus $G_{f_q}(z_0)=g_q(G_s(z_0))
\in\mathbb R$, an obvious contradiction.

\begin{figure}[h,t]
\includegraphics[width=10cm,height=3.5cm]{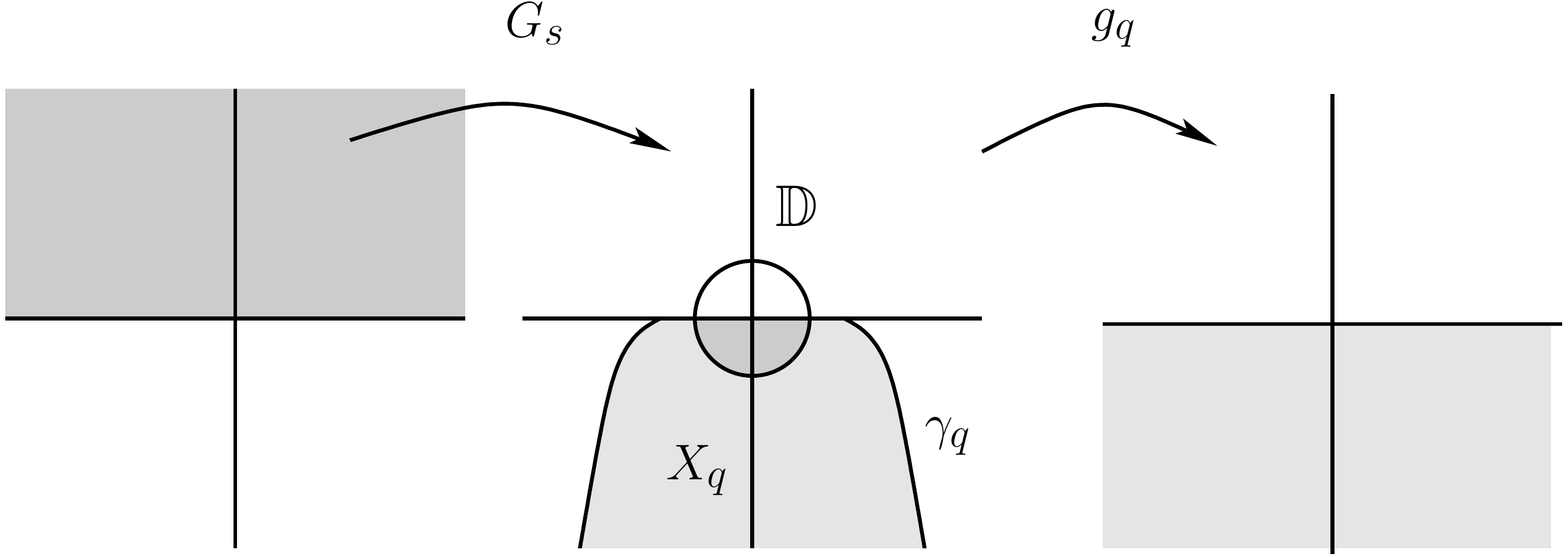}
\caption{The curve $\gamma_q$ is mapped by $g_q$ onto the real line.
The darker shadings correspond to each other via $G_s$, while $g_q$ maps
the entire shaded area (dark and light) onto the lower half-plane.}\label{fig1}
\end{figure}


Fix now $q\in(0,1)$.
It is trivial that $g_q$ has a unique inverse, call it $\Phi$, defined
around zero (more precise, at least on $M_q\mathbb D$), which fixes 
zero. We would like to 
extend this inverse to the whole lower half-plane; then it 
would easily follow that $\Phi(\mathbb C^-)=X_q$. The lack of finite
asymptotic values for $g_q$ guarantees that the only
impediment to such an analytic extension would be a finite
critical value of $g_q$. Thus, assume towards contradiction that
there exists a $c\in X_q$ (for precision, assume $\Re c>0$)
so that $g_q'(c)=0$. Without loss of generality
we may assume that $g_q(c)\in\mathbb C^-$
and $g_q(c)$ is the point closest to the origin in the lower half-plane
to which $\Phi$ does not extend analytically. We consider a 
half-line $r$ starting at $g_q(c)$ and going to infinity inside
$\mathbb C^-\cap i\mathbb C^-$ without cutting $g_q(c)\mathbb D$ and 
so that no $z\in r\setminus\{g_q(c)\}$ is a critical value for $g_q$
(possible since the set of critical values of $g_q$ is at most 
countable).
By the same Lemma \ref{curve}, we conclude that $r$ has a preimage 
$p_q$ via $g_q$ inside $X_q$, right of the imaginary axis, with
both ends at infinity, which determines a domain $D\subset X_q$ 
that does not contain $i(-\infty,0)$ and having boundary
$\partial D=p_q$. The choice for $p_q$ is made so that there is no 
other preimage in $X_q\setminus D$ of $r$ which is in the same connected component
of $g_q^{-1}(r)$ as $p_q$.

For the sake of clarity, we group most of the rest of the proof in
the following lemma. The reader will probably find Figure \ref{fig2} helpful in following it.

\begin{lemma}
\begin{enumerate}
\item There exists a point $b\in D$, $|b|>1$, so that $g_q(b)=0$.
\item For any $m>|c|+|g_q(c)|+1+|b|$, there exists $M\ge m$ and a path 
$\varpi$ in $\overline{D}$ uniting the two points of $p_q\cap
g_q^{-1}(\partial(M\mathbb D))$ so that $g_q(\varpi)\cap 
m\mathbb D=\varnothing$.
\item Let $\Pi_m$ be the path obtained by 
concatenating $\varpi$, the bounded part(s) of $p_q\setminus
(\varpi\cup\mathbb D)$ 
and, if existing, the segments $\partial\mathbb D\cap D$. Then $\Pi_m$
is a closed curve containing $b$ inside it and there exists $m>0$
so that $\{g_q(1/w)\colon w\in\Pi_m\}
\subset\mathbb C^+$ and $g_q(\Pi_m)\cap\{g_q(1/w)\colon w\in\Pi_m\}
=\varnothing$.
\item In particular, $w\mapsto g_q(w)$ and $w\mapsto g_q(w)-g_q(1/w)$
have the same number of zeros in the simply connected domain 
determined by the simple closed curve $\Pi_m$.
\end{enumerate}
\end{lemma}
\begin{proof}

Clearly, Lemma \ref{picard} guarantees that $g_q(D)\supseteq
\mathbb C\setminus r$, which contains zero. 
Thus, there exists a point $b\in D$ so that $g_q(b)=0$.
Moreover, let us recall that $G_{f_q}=g_q\circ G_s$
is the Cauchy transform of a probability measure, and, as such, it
maps the upper half-plane into the lower one. Thus,
as $G_s(\mathbb C^+\cup[-2,2])=\overline{\mathbb D}\cap\mathbb C^-$,
it follows that $g_q(w)\neq0$ for any $w\in
\overline{\mathbb D}\cap\mathbb C^-$, so $b\in\mathbb C^-
\setminus\overline{\mathbb D},$ so in particular $|b|>1$. This
proves (1).

To prove (2), recall that by Theorem \ref{DCA} the set
$\{z\in\mathbb C\colon|g_q(z)|>m\}$ is connected for any $m\ge0$ and
in particular for $m$ as in the statement of the lemma.
Also, from the construction of $r$ and $p_q$, it is clear that
$p_q\cap
g_q^{-1}(\partial(M\mathbb D))$ contains exactly two points for any
$M\ge m$. Assume now towards contradiction that there is no
path uniting those two points inside the set $D\setminus g^{-1}_q
(m{\mathbb D})$, and in particular, of course, in
$D\setminus g^{-1}_q
(m\overline{\mathbb D})$. Thus the open set 
$g^{-1}_q(m{\mathbb D})\cap D$ must contain an unbounded smooth path.
Choose such a path, and call it $\Gamma$. Choose $0<T=1+2\inf\{|z|
\colon z\in\Gamma\}$ (so that $\Gamma\cap T\mathbb D\neq\varnothing$)
and let $M>2\max\{|g_q(z)|\colon|z|\leq T\}+m$.
Since $g_q(\overline{z})=\overline{g_q(z)}$, it follows that 
$\{z\in\mathbb C\colon|g_q(z)|>M\}$ does not intersect the set
$T\mathbb D\cup\Gamma\cup\overline{\Gamma}$ (here again 
$\overline{\Gamma}$ denotes complex conjugate, not closure).
But this set disconnects $\mathbb C$: for example,
the set $p_q\cap\{z\in\mathbb C\colon|g_q(z)|>M\}$
contains two nonempty connected components separated by 
$T\mathbb D\cup\Gamma\cup\overline{\Gamma}$. This contradicts
Theorem \ref{DCA}. Thus a path $\varpi=\varpi_{m}$ as described
in our lemma must exist.

The fact that $\Pi_m$ exists and surrounds $b$ exactly once
is a trivial consequence of (1), (2) and the entireness of $g_q$.
Moreover, from $\Pi_m$'s construction, the set $\{1/w\colon w\in
\Pi_m\}\subset\mathbb C^+\cap\overline{\mathbb D}$. Thus, as noted
in the proofs of (1) and (2), $\{g_q(1/w)\colon w\in
\Pi_m\}\subset g_q(\mathbb D\cap\mathbb C^+)\subset\mathbb C^+$ (recall that $g_q(\overline{z})=
\overline{g_q(z)}$) is a bounded set for any $m$ (one can
choose the bound $\max g_q(\overline{\mathbb D})$, which is obviously 
independent of $m$ and depends only on $q$).
To finish the proof of (3) we only need to 
argue that for $m$ large enough, the set $g_q(\Pi_m)$ does
not intersect $\mathbb C^+\cap(\max g_q(\overline{\mathbb D})
\mathbb D)$. Indeed, if we have a point $z\in p_q$, then 
$g_q(z)\in r\subset\mathbb C^-$, while if $z\in\partial\mathbb D
\cap\mathbb C^-$, then again $g_q(z)\in\mathbb C^-.$
We only need to show that for $m$ large enough $\varpi=\varpi_m$
is mapped in the complement of $\max g_q(\overline{\mathbb D})
\mathbb D$. However, this follows from part (2) by simply choosing
$m\ge\max g_q(\overline{\mathbb D})$. This proves (3).

\begin{figure}[h,t]
\includegraphics[width=10cm,height=8cm]{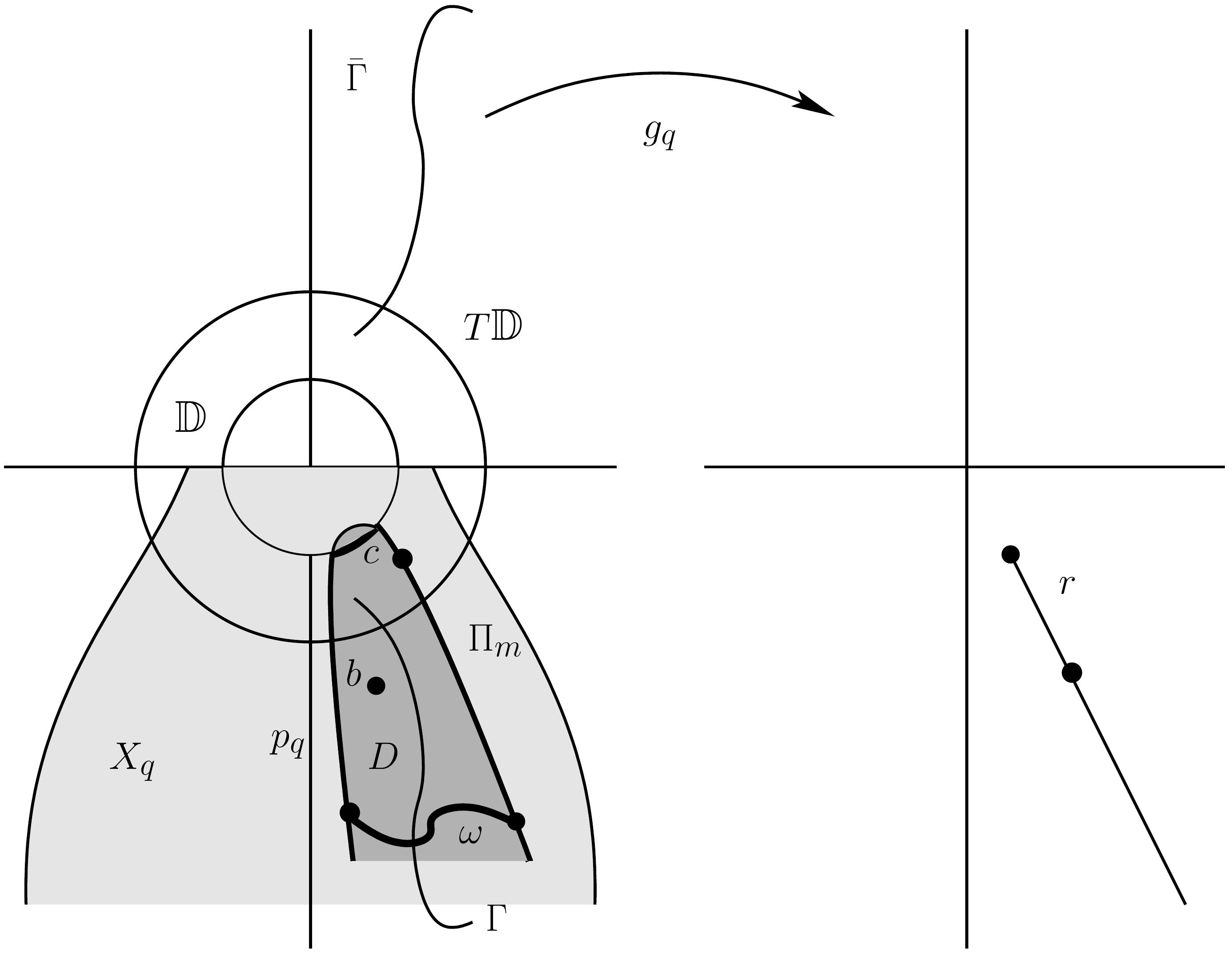}
\caption{The region shaded in dark grey is mapped by $g_q$ surjectively onto
$\mathbb C\setminus r$.}\label{fig2}
\end{figure}

In order to prove (4), we will apply Theorem \ref{var}
to $f(z)=-g_q(z)+g_q(1/z)$ and $g(z)=g_q(z)$.
The relation 
$|g_q(z)+(-g_q(z)+g_q(1/z))|<|g_q(z)-g_q(1/z)|+|g_q(z)|$
is equivalent to 
$\left|\frac{g_q(1/z)}{g_q(z)}\right|<\left|\frac{g_q(z)-g_q(1/z)}{
g_q(z)}\right|+1=\left|1-\frac{g_q(1/z)}{
g_q(z)}\right|+1,$ for $z\in\Pi_m$. Consider two cases: first, if
$z\in\varpi$, then $|g_q(z)|>m>|g_q(1/z)|$, so that 
$\left|\frac{g_q(1/z)}{g_q(z)}\right|-1<0<
\left|1-\frac{g_q(1/z)}{
g_q(z)}\right|$. Second case, if $z\in\Pi_m\setminus\varpi$,
then $g_q(z)\in\mathbb C^-$ and $g_q(1/z)\in\mathbb C^+$.
Thus, $g_q(z)$ and $g_q(1/z)$ cannot be positive multiples
of each other, i.e. $\frac{g_q(1/z)}{g_q(z)}\not\in[0,+\infty)$.
Generally, for the relation $|a|-1=|a-1|$ to hold it is necessary
that $a\ge1$. Applying this observation to $a=\frac{g_q(1/z)}{g_q(z)}$,
we conclude that the inequality must be strict also for $z\in\Pi_m
\setminus\varpi$.

Thus, we conclude that $g_q(z)$ and $g_q(z)-g_q(1/z)$ have the same
number of zeros in the domain delimited by $\Pi_m$. This 
completes the proof of (4) and of our lemma.
\end{proof}

The proof of our main theorem is now almost complete.
We will use equation $\eqref{real}$ to obtain a
contradiction. 
Our assumption that $g_q$ has a critical point $c$ in $X_q$ has led
us to conclude by part (1) of the previous lemma that the equation 
$g_q(z)=0$ has a solution $b\in D\subset X_q$. By part (4) of the
same lemma, the map
$D\ni w\mapsto k\Theta(w)=g_q(w)-g_q(1/w)$ must then have a zero 
in $D$. But $D\cap\mathbb R=\varnothing,$ and we have seen in Equation
\eqref{real} that the zeros of $\Theta$ are real. This is a 
contradiction.

We conclude that $g_q$ has no critical points in $X_q$, and so
$\Phi$ has an analytic continuation to the whole lower half-plane.
As noted at the beginning of the proof, this implies free
infinite divisibility for $f_q(x)dx$.
\end{proof}

\end{document}